\documentclass[reqno,12pt]{amsart}
\usepackage{amscd,amssymb}

\title{ On a Theorem of N. Katz and Bases in Irreducible Representations }
\author{ David Kazhdan}
\address{ Institute of Mathematics, The Hebrew University, Jerusalem, Israel}
\email{kazhdan@math.huji.ac.il}
\dedicatory{ Dedicated to the memory of Leon Ehrenpreis}
\begin{document}
\maketitle
\newcommand{\set}[1]{\left\{#1\right\}}

\newcommand{\un}{{\underline}}
\newcommand{\mR}{{\mathbb R}}
\newcommand{\mC}{{\mathbb C}}
\newcommand{\mZ}{{\mathbb Z}}
\newcommand{\mP}{\mathbb P}
\newcommand{\mG}{\mathbb G}
\newcommand{\mA}{\mathbb A}
\newcommand{\mF}{\mathbb F}
\newcommand{\mQ}{\mathbb Q}

\newcommand{\ho}{\hookrightarrow}

\newcommand{\Gg}{\gamma}
\newcommand{\GG}{\Gamma}
\newcommand{\bO}{\Omega}
\newcommand{\bo}{\omega}
\newcommand{\bl}{\lambda}
\newcommand{\bL}{\Lambda}
\newcommand{\bs}{\sigma}
\newcommand{\bS}{\Sigma}
\newcommand{\ep}{\epsilon}
\newcommand{\D}{\Delta}

\newcommand{\mcA}{\mathcal A}
\newcommand{\mcB}{\mathcal B}
\newcommand{\mcC}{\mathcal C}
\newcommand{\mcD}{\mathcal D}
\newcommand{\mcE}{\mathcal E}
\newcommand{\mcF}{\mathcal F}
\newcommand{\mcG}{\mathcal G}
\newcommand{\mcH}{\mathcal H}
\newcommand{\mcI}{\mathcal I}
\newcommand{\mcJ}{\mathcal J}
\newcommand{\mcK}{\mathcal K}
\newcommand{\mcL}{\mathcal L}
\newcommand{\mcM}{\mathcal M}
\newcommand{\mcN}{\mathcal N}
\newcommand{\mcO}{\mathcal O}
\newcommand{\mcP}{\mathcal P}
\newcommand{\mcQ}{\mathcal Q}
\newcommand{\mcR}{\mathcal R}
\newcommand{\mcS}{\mathcal S}
\newcommand{\mcT}{\mathcal T}
\newcommand{\mcU}{\mathcal U}
\newcommand{\mcV}{\mathcal V}
\newcommand{\mcW}{\mathcal W}
\newcommand{\mcX}{\mathcal X}
\newcommand{\mcY}{\mathcal Y}
\newcommand{\mcZ}{\mathcal Z}

\newcommand{\fa}{\mathfrak a}
\newcommand{\fb}{\mathfrak b}
\newcommand{\fc}{\mathfrak c}
\newcommand{\fd}{\mathfrak d}
\newcommand{\fe}{\mathfrak e}
\newcommand{\ff}{\mathfrak f}
\newcommand{\fg}{\mathfrak g}
\newcommand{\fh}{\mathfrak h}
\newcommand{\fii}{\mathfrak i}
\newcommand{\fj}{\mathfrak j}
\newcommand{\fk}{\mathfrak k}
\newcommand{\fl}{\mathfrak l}
\newcommand{\fm}{\mathfrak m}
\newcommand{\fn}{\mathfrak n}
\newcommand{\fo}{\mathfrak o}
\newcommand{\fp}{\mathfrak p}
\newcommand{\fq}{\mathfrak q}
\newcommand{\fr}{\mathfrak r}
\newcommand{\fs}{\mathfrak s}
\newcommand{\ft}{\mathfrak t}
\newcommand{\fu}{\mathfrak u}
\newcommand{\fv}{\mathfrak v}
\newcommand{\fw}{\mathfrak w}
\newcommand{\fx}{\mathfrak x}
\newcommand{\fy}{\mathfrak y}
\newcommand{\fz}{\mathfrak z}
\newcommand{\ti}{\tilde}
\newcommand{\va}{\vartriangleleft}
\newtheorem{theorem}{Theorem}[section]
\newtheorem{proposition}[theorem]{Proposition}
\newtheorem{conjecture}[theorem]{Conjecture}
\newtheorem{corollary}[theorem]{Corollary}
\newtheorem{lemma}[theorem]{Lemma}
\newtheorem{claim}[theorem]{Claim}
\newtheorem{maintheorem}[theorem]{Main Theorem}

\theoremstyle{definition}

\newtheorem{remark}[theorem]{Remark}
\newtheorem{example}[theorem]{Example}
\newtheorem{definition}[theorem]{Definition}
\newtheorem{algorithm}[theorem]{Algorithm}
\newtheorem{problem}[theorem]{Problem}
\newtheorem{pf}[theorem]{Proof}
\section{}

{Abstract}. N. Katz has shown that any irreducible representation of the Galois group of $\mF _q((t))$ has unique
extension to a {\it special} representation of the Galois group of $k(t)$  unramified outside $0$ and $\infty$
and tamely ramified at $\infty$. In this paper we analyze the number of not necessarily special such extensions
and relate this question to a description of bases in irreducible representations of multiplicative groups
of division algebras.\vspace{4mm}

Let $k=\mF _q,q=p^r$ be a finite field, $\bar k$ the algebraic closure of $k,F:=k((t))$ and $\bar F$ be
the algebraic closure of $F$. The restriction on $\bar k\subset \bar F$ defines a group homomorphism
$$Gal (\bar F/F)\to Gal (\bar k/k)=\hat \mZ$$
 and we define the {\it Weil group} of the field $F$ as the preimage
$\mcG _0\subset Gal (\bar F/F)$ of $\mZ \subset \hat \mZ$ under this homomorphism.\vspace{2mm}

 We denote by $\underline \mP ^1$ the projective line over $k$, set $E:=k(t)$ and denote by   $S$ the set of  points of
$\underline \mP ^1$. For any $s\in S$ we denote by $E_s$ the completion of $E$ at $s$. Using the parameter $t$ on
$\underline \mP ^1$ we identify the fields $E_0$ and $E_{\infty}$ with $F$ and therefore identify $\mcG _0$ with the  Weil groups of the fields   $E_0$ and $E_{\infty}$.\vspace{2mm}

 Let $\ti E$ be the maximal extension of the field $E$ unramified outside $0$ and $\infty$ and tamely ramified
at $\infty$.
We denote by  by $\mcG \subset Gal (\ti E/E)$ the Weil group corresponding to the  extension $\ti E/E$.  We have
the natural  imbeddings
$$\mcG _0\ho \mcG ,\mcG _{\infty}\ho \mcG$$
well defined up to  conjugation. Therefore for any complex representation $\rho$ of $\mcG$ the restrictions to $\mcG _0,\mcG _{\infty}$ define representations $\rho _0,\rho _{\infty}$ of the corresponding local groups.
The group $\mcG$ has a unique maximal quotient $\bar \mcG$ such that the Sylow $p$-subgroup of $\bar \mcG$ is normal. As shown by N.Katz  (\cite{Katz}) the
composition $\mcG _0\to \bar \mcG$ is an isomorphism.\vspace{2mm}

{\bf Remark} A finite-dimensional irreducible representation $\rho _0$ of $\mcG$ is called  {\it special} if it  factors through a representation of the group $\bar \mcG$.
One can restate the theorem of N.Katz by saying that for any irreducible representation $\rho _0$ of $\mcG _0$ there exists
a unique special representation $\rho _{sp}$ of the group $\mcG$ whose restriction to $\mcG _0$ is equivariant to $\rho _0$.\vspace{2mm}

 Let $D_0$ be a skew-field with  center $F,dim _FD_0=n^2$,  $G_0:=D_0^{*}$ be the multiplicative group of
$D_F$ and $\rho _0$  be  an  $n$-dimensional indecomposable  representation of the group $\mcG _0$.

\begin{definition}a) We denote by  $\ti \bs (\rho _0)$ the irreducible discrete series representation of the group $GL_n(F)$ which corresponds to $\rho _0$ under the local  Langlands correspondence ( see for example (\cite{H})
and by $\bs (\rho _0)$ the irreducible representation of the group $G_0$ which corresponds to $\ti \bs (\rho _0)$ as in \cite{DKV}.\vspace{2mm}

b) We  denote by $ r(\rho _0)$ the
formal dimension of the representation $\ti \bs (\rho _0)$ where the formal dimension is normalized in such a way that the formal dimension of the Steinberg representation is equal to $1$.
 Analogously for any  indecomposable representation $\rho _{\infty}$ of the group $\mcG _{\infty}$ we define an integer  $r (\rho _{\infty})$.\vspace{2mm}

c) We denote by  $A(\rho _0)$ the set of equivalence classes of  $n$-dimensional irreducible representations
 $\rho$ of the group $\mcG$ whose restriction to $\mcG _0$ is equivalent to $\rho _0$ and the  restriction to $\mcG _{\infty}$ is indecomposable.\vspace{2mm}

\end{definition}
\begin{theorem} For any $n$-dimensional irreducible $\bar \mQ _l$-representation  of the group $\mcG _0$ the sum $\sum _{\rho \in  A  (\rho _0)}r(\rho _{\infty})$  is equal to $ r(\rho _0)$.
\end{theorem}
{\bf Proof.} Let $\mA=\prod _{s\in S}E_s$ the ring of adeles of $E$ and $D$ be a  skew-field  with  center $E$
unramified outside $\{ 0,{\infty}\},D_0:=D\otimes _E E_0$ and $D_{\infty}:=D\otimes _E E_{\infty}$.
Then $D_0,D_{\infty}$ are local skew-fields. Let $\underline G$ be  the multiplicative group of $D$
considered as an the algebraic $E$-group. \vspace{2mm}

It  follows from \cite{L} that we can identify the set
$A(\rho _0)$ with the set of automorphic representations $\ti \pi =\prod _{s\in S}\ti \pi _s$ of the group $GL_n(\mA )$ such that the representation
$\ti {\pi _0}$ is equivalent to $\ti \bs (\rho _0)$ and the representation $\ti \pi _{\infty}$ is of discrete series. Then it follows from \cite{DKV} that  we  can identify the set  $A(\rho _0)$ with the set of automorphic representations $\pi =\prod _{s\in S}\pi _s$  of the group $\underline G(\mA )$ such that the representation
$\pi _0$ is equivalent to $\bs (\rho _0)$.  We will use this identification for the proof of the Theorem 1.2.\vspace{2mm}

 We denote by  $N:D_0\to F$  the reduced norm and define
$$\mu :=\nu \circ N:D_0^{*}\to \mZ ,K_0:=\mu ^{-1}(0)$$

 where $\nu :F^{*}\to \mZ$ is the valuation.
Then $K_0\subset D_0^{*}$ is a maximal compact subgroup. We define the
{\it first congruence subgroup} $K^1_0$ by
$$K^1_0:=\{k\in K_0|\mu (k-Id)>0\}$$
 As is well known $K^1_0$
 is a normal subgroup of $D_0^{*}$ such that $K_0/K^1_0=\mF _{q^n}^{*}$ and
 $D_0^{*}/K^1_0=\mZ \ltimes \mF _{q^n}^{*}$ where $\mZ$ acts on $\mF _{q^n}^{*}$ by
$(n,x)\to x^{q^n}$. \vspace{2mm}

For any $s\in S-\{0,{\infty} \}$ we identify the group $G_{E_s}$ with $GL(n,E_s)$  and define $K_s:=GL(n,\mcO _s)$.
We  write
$G_\mA :=D_\infty ^{*}\times GL_n (\mA ^0)$ where
$$GL_n (\mA ^\infty ):= D_0 ^{*}\times \prod _{s\in S-\{0 ,{\infty} \}}GL(n,E_s)$$
 and define
$$K^0:=\prod _{s\in S-\{0 ,{\infty} \}}K_s\times K_{E_{\infty}},
K^1:=\prod _{s\in S-\{0 ,{\infty} \}}K_s\times K^1_{E_{\infty}}$$
where   $K^1_{E_{\infty}} \subset K_{E_{\infty}} \subset D_{\infty} ^{*}$ is the first congruence subgroup of $G_{\infty}$. \vspace{2mm}

For any  irreducible representation  $\pi$ of the group $G_0$   we denote by $\ti \pi$ the discrete series representation of the group $GL_n(F)$ corresponding to $\ti \pi$ as in \cite{DKV}.

\begin{lemma} a)  For any irreducible complex representation
$\kappa :D_0^{*}/K^1_0\to
Aut(W)$ and any character $\chi :K_0/K^1_0\to \mC ^{*}$ we have
$$dim (W^\chi )\leq 1$$
 where $W^\chi =\{ w\in W| \kappa (k)w=\chi (k)w,k\in  K_0\}$.\vspace{2mm}

b) For any  irreducible representation  $\pi$ of the group $G_0$  the formal dimension of $\ti \pi$ is equal to the dimension of $\pi$.
\end{lemma}
{\bf Proof}. Part a) follows from the isomorphism
$D_0^{*}/K^1_0=\mZ \ltimes \mF _{q^n}^{*}.$\vspace{2mm}

Part b) follows from \cite{DKV}$.\square$\vspace{2mm}

We see that the following equality implies the validity of the Theorem 1.2.
\begin{claim}  For any $n$-dimensional irreducible $\bar \mQ _l$-representation of the group $\mcG _0$ the sum $\sum _{\pi \in  A (\rho _0)}dim (\pi  _{\infty})$  is equal to $dim(\bs (\rho _0))$.
\end{claim}
The proof of  Claim is based on the following result.

\begin{proposition} The  product map $D_0^{*}\times K^1\times G_E\to G_\mA$
 is a bijection.
\end{proposition}
{\bf Proof of the Proposition}. The surjectivity follows from Lemma 7.4 in \cite{HK}. To show the injectivity
it is sufficient to check the equality
$$(D_0^{*}\times K^1)\cap G_E=\{ e\}$$
 which is obvious.$\square$ \vspace{2mm}

We denote by $\mC (G_\mA /G_E)$ the space of
locally constant functions on $G_\mA /G_E$ with compact support, by $\mC (G_0)$ the space of
locally constant functions on $G_0$ with compact support and by $L\subset  \mC (G_\mA /G_E)$ the subspace of
$K^1$-invariant functions. The group  $G_0\times D^\star _\infty /K^1_\infty $ acts naturally on $L$.\vspace{2mm}

Let $\rho _0$ be  an indecomposable representation of the group $\mcG _0$. We denote by  $(\bs (\rho _0),V(\rho _0))$ the corresponding representation of the group $G_0$ and
identify the set $ A(\rho _0)$ with  the set of automorphic representations $\pi ^a=\prod _{s\in S}\pi ^a_s$  of the group $\underline G(\mA )$ such that the representation
$\pi ^a_0$ is equivalent to $\bs (\rho _0)$ and  the representation
$\pi ^a_\infty$ is trivial on $K^1_\infty$. Let
$$\mcH: =\prod _{s\in S-\{0 ,{\infty} \}}\mcH _s$$
 where $\mcH _s$ is the spherical Hecke algebra for $G(F_s)=GL(n,F_s)$. By construction, the commutative algebra $\mcH$ acts on the
$D^\star _0\times D^\star _{\infty} /K^1$-module $L$. For any $a\in A(\rho _0)$ we define
$$L_a:=Hom_{G_\mA ^\infty} (\pi ^a, \mC (G_\mA /G_E))=Hom_{G_0\times \mcH }(\bs (\rho _0),L)\subset Hom_{G_0 }(\bs (\rho _0),L)$$
\begin{lemma}a) The restriction $r:L\to \mC (G_0)$ is an isomorphism of $G_0$-modules where
 $G_0$ acts on $\mC (D_0^{*})$ by  left translation.\vspace{2mm}

b) $Hom_{G_0 }(\bs (\rho _0),L)=V^\vee$ where $V^\vee$ is the dual space to $V(\rho _0)$.\vspace{2mm}

c) $V^\vee =\oplus L_a,a\in A(\rho _0)$ where the algebra  $\mcH$ acts on $L_a,a\in A(\rho _0)$ by a character $\chi _a:\mcH \to \bar \mQ _l^\star ,\chi _a\neq \chi _{a'}$ for
$a\neq a'$ and the representations $\pi _\infty^a$ of the group $D^\star _{\infty} /K^1$ on $M_a$ are irreducible.\vspace{2mm}

d) The representations $\pi _\infty^a$ are  associated with the restriction $\rho (a)_{\infty}$ by the local Langlands correspondence.\vspace{2mm}
\end{lemma}

{\bf Proof}. The Lemma follows immediately from the Proposition and the strong multiplicity one theorem (\cite{PS}and \cite{DKV}).$\square$\vspace{2mm}

This Lemma  implies the validity of  Claim and therefore of  Theorem 1.2.
Indeed we have

$$dim (V)=dim (V^\vee )=\sum _{a\in A(\rho _0)}dim (L_a)=\sum _{a\in A(\rho _0)}dim (\pi _\infty^a)=
\sum _{a\in A(\rho _0)} r(\rho (a)_{\infty} )\square$$\vspace{4mm}

One can ask whether one can extend Theorem 1.2 to the case of other groups. More precisely,
let $G$ be a split reductive  group with a connected center and ${}^LG$ be the Langlands dual group. Consider a homomorphism $\rho _0:\mcG _0\to {}^LG$  such that the connected component of
the centralizer $Z_\rho :=Z_{{}^LG}(Im(\rho ))$ is unipotent. Let $[Z_\rho]$ be the group of
connected components of the centralizer $Z_\rho$.  Conjecturally, one can associate with  $\rho _0$ an $L$-packet of irreducible representations $\pi _\rho (\tau )$ of the group $G_0:=G(F)$ parameterized by
irreducible representations $\tau$ of $[Z_\rho]$ and there exists an integer $r(\rho _0)$ such that the formal dimension of  $\pi _\rho (\tau )$ is equal to $r(\rho )dim (\tau)$.\vspace{2mm}

We denote by $A^G(\rho _0)$ the set of conjugacy  classes of homomorphisms $\rho :\mcG \to {}^LG$ whose restriction on $\mcG _0$ is conjugate to $\rho _0$ and such that the connected component of
the centralizer of the restriction on $\mcG _\infty$ is unipotent.\vspace{2mm}

{\bf Question}. Is it true that $r(\rho _0)=\sum _{a\in A(\rho _0)}r(\rho _{\infty})$ where $r(\rho _{\infty})$ is defined in the same way as $r(\rho _0)$?\vspace{2mm}

\section{}
Let $G$ be a  reductive group over a local field. As is well known one can realize the spherical Hecke algebra $\mcH$
of $G$ {\it geometrically}, that is  as the Grothendick group of the monoidal category of perverse sheaves
on the affine Grassmanian. Analogously in the case when  $G$ be a  reductive group over a global field of
positive characteristic the   unramified geometric Langlands conjecture predicts the existence of a
 geometric realization of the corresponding space of automorphic functions.\vspace{2mm}

Let $\underline C$ be a smooth absolutely irreducible $\mF _q$-curve, $q=p^m$,
$S$ be the set of geometric points  of
$\underline C$, $\GG :=\pi _1 (\underline C)$. For any $s\in S$ we denote by $Fr_s \subset \GG$ the conjugacy
class of the Frobenius at $s$.\vspace{2mm}

Let  $E$ be the field of rational functions on $\underline C$.
For any $s\in S$ we denote by $E_s$ the completion of $E$ at $s$ and we denote by
$\mA$ be the  ring of
adeles of $E$. Fix a prime number $l\neq p$.\vspace{2mm}

Let $\underline G$ be a split reductive group,
and $\hat K :=\prod _{s\in S}G( \mcO _s)\subset G(\mA)$ be the standard maximal
compact subgroup. An irreducible representation
$(\pi ,V) =\otimes _{s\in S}(\pi _s ,V_s)$ of $G(\mA )$ is {\it unramified}
if $V^{\hat K}\neq \{ 0\}$. In this case $dim (V^{\hat K})=1$. So for any
unramified  representation $(\pi ,V)$ of the group $ G(\mA)$ there is a special
{\it spherical } vector $v_{sp}\in V$ defined up to a multiplication by a scalar. \vspace{2mm}

Let ${}^LG$ be the Langlands dual group and
$\rho$  a homomorphism from $\GG$ to ${}^LG (\bar \mQ _l)$ such that for any $s\in S$ the conjugacy
class $\Gg _s:=\rho (Fr_s)\subset {}^LG (\bar \mQ _l)$ is semisimple. In such a case we can define
unramified representations $(\pi _{\Gg _s}, V_s)$ of local groups $G(E_s)$ and the  representation
  $(\pi (\rho ),V_\rho )=\otimes _s (\pi _{\Gg _s}, V_s)$ of the  adelic group $ G(\mA)$. According
to the  unramified geometric Langlands conjecture the homomorphism $\rho$  defines [at least in the case when
$\rho$ is tempered]  an imbedding
$$i_\rho :V_\rho \to \bar \mQ _l (K\backslash G(\mA)/G(E))$$
and  a function $f_\rho :=i_\rho (v_{sp})$ which  is defined up to a multiplication by a scalar.
\vspace{2mm}

We can identify the set $K\backslash G(\mA)/G(E)$ with the set of $\mF _q$-points of
the stack $\mcB _G$  of principal $G$-bundles on $\underline C$ and   the  unramified geometric Langlands
correspondence  predicts the existence of a perverse Weil  sheaf $\mcF (\rho )$
on  $\mcB _G$ such that the function  $f_\rho $ is given by  the trace of the
Frobenius automorphisms on stalks of $\mcF (\rho )$. ( See  \cite{G})\vspace{2mm}

If one considers ramified  automorphic representations
$(\pi ,V) =\otimes _{s\in S}(\pi _s ,V_s)$ of $G(\mA )$ then there is no
natural way to choose a special vector in $V$. So on the "geometric" side one expects not an object
$\mcF (\rho )$
but an abelian  category $\mcC (\rho )$ which is a product of local categories $\mcC (\rho _s)$
such that the Grothendick K-group of the category  $[\mcC (\rho _s)]$ coincides with
the subspace $V^0_s$ of the  minimal $K$-type vectors  of the space $V_s$ of the local representation. Such
geometric realization of the space $V^0_s$ would define a special basis of vector spaces $V^0_s$
which  would be a non-archimedian analog of  Lusztig's
canonical basis. Here we consider only
the case of an anisotropic group when the  minimal $K$-type subspace  $V^0_s$ coincides with the space $V_s$
of the representation of $G$. Moreover we will only discuss a slightly weaker data of  a {\it projective basis}
where a projective basis in a finite-dimensional vector space $T$ is
a  decomposition of the space $T$  in a direct sum of one-dimensional subspaces.
So one could look  for a special basis of vector spaces $V_s$
which  would be a non-archimedian analog of the Lusztig's
canonical basis. \vspace{2mm}

Let as before $F:=k((t)),D_0$ be a skew-field with  center $F,dim _0D_0=n^2,G_0$ be
the multiplicative group of  $D_0$ and $\bs :G_0\to Aut (V)$  a complex
irreducible   continuous representation of the group  $G_0$.

\begin{theorem} For any irreducible representation $\tau :D_F^{*}\to Aut (T)$  of the group  $D_F^{*}$
there exists  a "natural" projective basis $=\oplus _aT_a$ of $T$.
\end{theorem}
\begin{remark} The construction is global. In particular I don't know how to
define  a  projective basis in the case when  $F$
is  a local  field of characteristic zero. It would be very interesting to find
a local construction of a  projective basis.
\end{remark}
{\bf The construction.} As follows from Lemma 1.6 c) we have a decomposition $V^\vee =\sum _{a\in  A(\rho _0)}M_a$
where the group $D^\star _{\infty} /K^1_{\infty}$ acts irreducibly on $M_a$. Therefore the group $\mF _{q^n}^\star =K_{\infty}/K^1_{\infty}$
acts on $M_a$ and we have a decomposition of $M_a$ into the sum of eigenspaces for the action of the group $\mF _{q^n}^\star$.
As follows from Lemma 3 a) these eigenspaces  are one-dimensional.

\end{document}